\documentclass{amsart}
\usepackage[T1]{fontenc}
\usepackage[latin1]{inputenc}
\makeatletter

\theoremstyle{plain}    
\newtheorem{thm}{Theorem}[section]
\numberwithin{equation}{section} 
\numberwithin{figure}{section} 
\theoremstyle{plain}    
\newtheorem{lem}[thm]{Lemma} 
\theoremstyle{plain}    
\theoremstyle{plain}    
\newtheorem{Def}[thm]{Definition} 
\theoremstyle{remark}

\theoremstyle{remark}

\usepackage[arrow,matrix,curve,ps]{xy}
\xyoption{dvips}
\CompileMatrices

\makeatother

\begin{document}

\title{A Weak Kawamata-Viehweg Vanishing Theorem on compact K\"ahler 
       manifolds}


\author{Thomas Eckl}



\thanks{}

\address{Thomas Eckl, Institut für Mathematik, Universität
  Bayreuth, 95440 Bayreuth, Germany}

\email{thomas.eckl@uni-bayreuth.de}

\urladdr{http://btm8x5.mat.uni-bayreuth.de/\~{}eckl}

\maketitle

\pagestyle{myheadings}
\markboth{THOMAS ECKL}{A WEAK KAWAMATA-VIEHWEG VANISHING}

\bibliographystyle{alpha}

\section{Introduction}

\noindent
Using a subtle version of the Bochner technique, J.-P. Demailly and Th. 
Peternell were able to prove one instance of the Kawamata-Viehweg Vanishing
Theorem for (even singular) K\"ahler manifolds:
\begin{thm}[{\cite[Thm. 0.1]{DP02}}]
Let
$X$
be a normal compact K\"ahler space of dimension
$n$
and
$L$
a nef line bundle on
$X$. 
Assume that
$L^2 \neq 0$.
Then
\[ H^q(X, K_X + L) = 0 \]
for 
$q \geq n-1$.
\end{thm}

\noindent
The aim of this paper is to use the same methods to prove a weak version of
Kawamata-Viehweg Vanishing on compact K\"ahler manifolds for all
$q > n - \nu(L)$, 
which also works for
pseudo-effective line bundles. It is only a weak version, because we have to 
tensorize 
$K_X + L$
with the \textbf{upper regularized multiplier ideal sheaf} of a singular 
hermitian metric on
$L$:
\begin{Def}
Let
$X$
be a compact K\"ahler manifold of dimension
$n$
and
$L$
a pseudo-effective line bundle on
$X$.
Let
$h_{\min}$
be a hermitian metric with minimal singularities among all positive singular 
hermitian metrics on
$L$.
Then the \textbf{upper regularized multiplier ideal sheaf}
$\mathcal{J}_+(L)$
is defined as
\[ \mathcal{J}_+(L) := \bigcup_{\epsilon \rightarrow 0} 
                                        \mathcal{J}(h_{\min}^{1+\epsilon}).\]
\end{Def}

\noindent
This multiplier ideal 
$\mathcal{J}_+(L)$
is certainly not optimal: there are examples of nef line
bundles where it is not trivial, see \cite[Ex.1.7]{DPS94}. At least, it is
conjectured (and true in dimension 1 and 2) that it equals the ordinary 
multiplier ideal
$\mathcal{J}(L) := \mathcal{J}(h_{\min})$.

\noindent
So the main result of this paper will be
\begin{thm} \label{KV+-thm}
Let
$X$
be a compact K\"ahler manifold of dimension
$n$
and
$L$
a pseudo-effective line bundle on
$X$
of numerical dimension
$\nu = \nu(L)$.
Then
\[ H^q(X, \mathcal{O}(K_X + L) \otimes \mathcal{J}_+(L)) = 0 \]
for
$q \geq n+1-\nu(L)$.
\end{thm}

\noindent
It will be proven as a corollary of
\begin{thm} \label{KV0-thm}
Let
$X$
be a compact K\"ahler manifold of dimension
$n$
and
$L$
a pseudo-effective line bundle on
$X$
of numerical dimension
$\nu = \nu(L)$
with a positive hermitian metric 
$h_{\min}$
with minimal singularities on
$L$.
For every
$\epsilon^\prime > 0$
there exists an
$0 < \epsilon < \epsilon^\prime$
such that the homomorphism
\[  H^q(X, \mathcal{O}(K_X + L) \otimes 
    \mathcal{J}(h_{\min}^{1+\epsilon^\prime}))
    \rightarrow H^q(X, \mathcal{O}(K_X + L) \otimes 
                                 \mathcal{J}(h_{\min}^{1+\epsilon})) \]
induced by the inclusion
$\mathcal{J}(h_{\min}^{1+\epsilon^\prime}) \subset 
 \mathcal{J}(h_{\min}^{1+\epsilon})$
vanishes for 
$q \geq n+1-\nu(L)$.
\end{thm}

\noindent
This theorem implies theorem~\ref{KV+-thm}, since the ascending chain of ideal
sheaves 
$\mathcal{J}(h_{\min}^{1+\epsilon})$,
$\epsilon \rightarrow 0$,
gets stable at some point by the Noetherian property, hence for
$\epsilon^\prime$
small enough
\[ \mathcal{J}(h_{\min}^{1+\epsilon^\prime}) = 
   \mathcal{J}(h_{\min}^{1+\epsilon}) = \mathcal{J}_+(L). \]
But then the homomorphism between the cohomology groups becomes an 
isomorphism, and the involved vector spaces must be 
$0$.

\noindent
The Bochner technique consists of using the Bochner-Kodaira-Nakano inequality
to prove the vanishing of certain cohomolgical classes, see \cite[\S 4]{Dem00}
for an application to classical vanishing theorems. The main obstacle to use 
this technique in our case is that the Bochner inequality is only true for 
\textit{smooth} metrics. This was circumvented by Demailly in \cite{Dem82}
who observed that for a compact K\"ahler manifold
$X$
and
$Z \subset X$
an analytic subset, 
$X \setminus Z$
has a complete K\"ahler metric. Using H\"ormander's confirmation of the Bochner
inequality for complete metrics (and elements of certain function spaces) in
\cite{Hoer65} Demailly constructed sequences of complete metrics converging to
the original K\"ahler metric and got vanishing results by going to the limit.
This limit process is successful despite of the metrics changing all the time,
since there is a uniform bound for all the occuring norms, and because one can
compare the different metrics.

\noindent
So in our case we start with a careful construction of singular hermitian 
metrics
$h_\epsilon$
on
$L$
which are smooth outside an analytic subset 
$Z_\epsilon$
(section~\ref{metric-sec}). They are composed of a (sufficiently small) part
controlling the eigenvalues of the curvature form (they are essential for the
Bochner inequality) and another part controlling the multiplier ideal sheaf 
(and hence the singularities) of the metric. The first is produced by applying
the Calabi-Yau theorem as in Boucksom's thesis \cite{Bou02}, the latter is 
constructed with the equisingular approximation theorem of \cite{DPS00}. 

\noindent
Then we use the Bochner technique for complete 
metrics on
$X \setminus Z_\epsilon$
converging to the starting K\"ahler metric
$\omega$
and go to the limit, using the uniform estimate in section~\ref{uniesti-sec}
and comparing the different metrics following the results in 
section~\ref{compmet-sec} (which more or less repeat the inequalities in 
\cite[\S 3]{Dem82}).

\vspace{0.2cm}

\noindent
\textit{Remark}. 
Since the technical details of the strategy described above are quite 
intricate and treated in a very summary way in \cite{DP02}, the author decided
to give all 
the steps in full details, for his own safety and for the convenience of the 
not so experienced reader -- of course without claiming any originality. The 
expert
may skip section~\ref{compmet-sec} alltogether and skim over the functional 
analytic details in section~\ref{Bochner-sec}.

\section{The construction of the metrics}
\label{metric-sec}

\noindent
As explained in the introduction, we construct metrics
$\hat{h}_\epsilon$
for arbitrarily small 
$\epsilon$
composed of two parts, and the first part is produced by using
\begin{thm}[Approximative Singular Calabi-Yau theorem] \label{AnSingCY-prop}
Let
$X$
be a compact K\"ahler manifold of dimension
$n$
with K\"ahler form
$\omega$
such that
$\int_X \omega^n = 1$,
and let  
$\alpha \in H^{1,1}(X,\mathbb{R})$
be a big class. Then for every
$\epsilon > 0$
there exists a closed positive current
$T_{\epsilon} \in \alpha$
with analytic singularities such that
\[ T_\epsilon(x)^n \geq (1-\delta_\epsilon)v(\alpha)\omega(x)^n \]
almost everywhere,
$\delta_\epsilon \rightarrow 0$
if
$\epsilon$
tends to
$0$, 
and the multiplier ideal
$\mathcal{J}(T_\epsilon)$
contains
$\mathcal{J}(T_{\min})$
for a positive current 
$T_{\min}$
in
$\alpha$
with minimal singularities. 
\end{thm}

\noindent
Here, 
$v(\alpha)$
denotes the volume of the class
$\alpha$,
defined by Boucksom \cite[3.1.6]{Bou02} as
\[ v(\alpha) := \sup_T \int_X T_{ac}^n, \]
where the
$T$'s
run through all closed positive currents in
$\alpha$
and
$T_{ac}$
is the absolute continuous part of the Lebesgue decomposition of
$T$
(\cite[3.1.1]{Bou02}).
Since
$\alpha$
is big, there is a closed positive current in
$\alpha$
bigger than 
$\epsilon\omega$
for some 
$\epsilon$,
and
$v(\alpha)$
is positive. 

\noindent
The proof of the theorem is contained in Boucksom's construction of
the current
$T$
(with arbitrary singularities) solving the Monge-Amp\`ere equation
\[ T_{ac}(x)^n = v(\alpha)\omega(x)^n \]
almost everywhere (\cite[Thm. 3.1.23]{Bou02}). His idea was of course to use 
the ordinary form of the Calabi-Yau theorem (where 
$\alpha$
contains a K\"ahler \textit{form} and 
$T$
will be a form satisfying the Monge-Amp\`ere equation 
\textit{everywhere}). To be able to do this, he proved a singular version of 
Fujita's theorem about the approximative Zariski 
decomposition (\cite[Thm. 3.1.24]{Bou02}):
\begin{thm}[Singular Fujita decomposition] \label{SingFujDec-thm}
Let
$\alpha \in H^{1,1}(X,\mathbb{R})$
be a big class. For all
$\epsilon > 0$
there exists a sequence of blow ups with smooth centers
$\mu: \widetilde{X} \rightarrow X$
and a decomposition
\[ \mu^\ast \alpha = \beta + \{E\}, \]
where
$\beta$
is a K\"ahler class and
$E$
is an effective
$\mathbb{R}-\!$
divisor such that
$|v(\alpha) - v(\beta)| < \epsilon$
and
$\mathcal{J}(\mu_{\ast} [E_k]) \supset \mathcal{J}(T_{\min})$.
\end{thm}
\begin{proof}
The additional statement compared to \cite[Thm. 3.1.24]{Bou02} is the inclusion
of the multiplier ideal sheaves. To get it, we must go through the proof of 
Boucksom: It starts with a sequence of K\"ahler currents
$T_k \in \alpha$
such that
\[ \int_X T^n_{k,ac} \rightarrow v(\alpha). \]
These
$T_k$
may be chosen in such a way that
$\mathcal{J}(T_k) \supset \mathcal{J}(T_{\min})$
(see for example the proof of Theorem 3.16 in \cite{Eck03}: The currents 
replacing the original approximating currents have potentials which are not 
less than the original ones).

\noindent
Next, we resolve the singularities of
$T_k$
and find a sequence of blow ups in smooth centers,
$\mu_k: \widetilde{X} \rightarrow X$,
such that the Siu decomposition (\cite[(2.18)]{Dem00})
\[ \mu^\ast T_k = [E_k] + R_k \]  
consists of the integration current of an effective
$\mathbb{R}$-\!
divisor 
$E_k$
and a smooth positive residue current
$R_k$. 
Since
$T_k$
is big, 
$\mu^\ast T_k$
is also big \cite[Prop.1.2.5]{Bou02}, hence there is a
$\delta > 0$
such that
\[ \mu^\ast T_k = R_k + [E_k] \geq \delta\omega. \]
This inequality remains valid, when we subtract an integration current of 
divisors, hence
$R_k \geq \delta\omega$.
But
$R_k$
is smooth, hence is a K\"ahler form. Furthermore, we have
$\int_X T_{k,ac}^n = \int_{\widetilde{X}} R_k^n$.
Finally, as 
$T_k$
is a current with analytic singularities, its potential
may be locally written as
\[ \phi_k = \theta_k + c \cdot \log(\sum|f_i|^2), \]
where
$\theta_k$
is a
$\mathcal{C}^\infty$
function, and an easy calculation shows that
\[ \mu^\ast \phi_k = \mu^\ast \theta_k + c \cdot \log(\sum|f_i^\prime|^2) + 
   c \cdot \log |g_k|^2, \]
where 
$g_k$
is a local equation of
$E_k$
and
$\sum|f_i^\prime|^2$
vanishes nowhere. Consequently,
$\mathcal{J}(\mu_{\ast} [E_k]) = \mathcal{J}(T_k)$,
and 
\[ \mu^\ast \alpha = \{R_k\} + \{E_k\} \]
is the desired decomposition, if we choose
$k$
big enough.
\end{proof}

\noindent
To prove the approximative singular Calabi-Yau theorem, let
us take a modification
$\mu: \widetilde{X} \rightarrow X$
and a decomposition
$\mu^\ast \alpha = \beta + \{E\}$
belonging to
$\epsilon$
as above. Let
$\widetilde{\omega}$
be a K\"ahler form on
$\widetilde{X}$, 
and set
$\widetilde{\omega}_\delta = \mu^\ast \omega + \delta \widetilde{\omega}$.
For all
$\delta > 0$,
$\widetilde{\omega}_\delta$
is a K\"ahler form on
$\widetilde{X}$, 
hence by the usual Calabi-Yau theorem, we can find a K\"ahler form 
$\theta_{\epsilon,\delta} \in \beta$
such that
\[ \theta_{\epsilon,\delta}(x)^n = 
   \frac{v(\beta)}{\int \widetilde{\omega}_\delta^n} \widetilde{\omega}(x)^n \]
for all 
$x \in \widetilde{X}$.
Now set
$T_\epsilon := \mu_\ast(\theta_{\epsilon,\delta} + [E_\epsilon])$ 
which is a closed positive current with analytic singularities in
$E_\epsilon$.
Furthermore, choosing
$\delta$
small enough and using the properties of
$\beta$
and
$\theta_{\epsilon,\delta}$,
we see that
\[ T_\epsilon(x)^n \geq (1-\delta_\epsilon)v(\alpha)\omega(x)^n \]
almost everywhere, and
$\delta_\epsilon \rightarrow 0$
if
$\epsilon$
tends to
$0$.
Since the multiplier ideals only depend on 
$[E_\epsilon]$,
the inclusion of multiplier ideals in the theorem remains true. 
\hfill $\Box$

\vspace{0.2cm}

\noindent
The construction of the second part of
$\hat{h}_\epsilon$
uses
\begin{thm}[Equisingular Approximation]
Let
$T = \alpha + i\partial\overline{\partial}\phi$
be a closed 
$(1,1)$-\!
current on a compact hermitian manifold
$(X,\omega)$,
where
$\alpha$
is a smooth closed
$(1,1)$-\!
form and
$\phi$
a quasi-plurisubharmonic function. Let
$\gamma$
be a smooth real 
$(1,1)$-\!
form such that
$T \geq \gamma$. 
Then one can write
$\phi = \lim_{\nu \rightarrow +\infty} \phi_\nu$
where
\begin{itemize}
\item[(a)]
$\phi_\nu$
is smooth in the complement 
$X \setminus Z_\nu$
of an analytic subset
$Z_\nu \subset X$;
\item[(b)]
$(\phi_\nu)$
is a decreasing sequence, and
$Z_\nu \subset Z_{\nu+1}$
for all
$\nu$;
\item[(c)]
for every
$t > 0$
\[ \int_X (e^{-2t\phi} - e^{-2t\phi_\nu})dV_\omega \]
is finite for
$\nu$
large enough and converges to
$0$
as
$\nu \rightarrow +\infty$;
\item[(d)]
$\mathcal{J}(t\phi_v) = \mathcal{J}(t\phi)$
for
$\nu$
large enough (``equisingularity'');
\item[(e)]
$T_\nu = \alpha + i\partial\overline{\partial}\phi_\nu$
satisfies
$T_\nu \geq \gamma - \epsilon_\nu \omega$,
where
$\lim_{\nu \rightarrow +\infty} \epsilon_\nu = 0$.
\end{itemize}
\end{thm}
\begin{proof}
See \cite{DPS00} or \cite[(15.2.1)]{Dem00} and especially the remark after the
proof. 
\end{proof}

\noindent
Now, let
$X$
be a compact n-dimensional K\"ahler manifold with K\"ahler form 
$\omega$
and let
$L$
be a holomorphic line bundle having a hermitian metric
$h_\infty$
with a curvature form 
$\Theta_{h_\infty}(L)$
of arbitrary sign. Let
$T_\epsilon$
be a closed current with analytic singularities in
$c_1(L)[-\epsilon\omega]$
such that
\[ (T_\epsilon + \epsilon\omega)^n \geq 
   \frac{v(c_1(L) + \epsilon\omega)}{2}\omega^n, \]
as in the approximative singular Calabi-Yau theorem above. There exists a 
hermitian metric
$h_\epsilon = h_\infty e^{-2\phi_\epsilon}$
on
$L$, 
such that
\[ \Theta_{h_\epsilon}(L) = T_\epsilon \]
(\cite{Bon95},\cite[Lem.4.1]{Eck03}). Next, we observe that
\[ v(c_1(L)+\epsilon\omega) \geq \epsilon^{n-l} \cdot 
                             (c_1(L)^l.\omega^{n-l})_\geq 0 \]
for all
$0 \leq l \leq n$
(\cite[p.86]{Bou02}), hence there is a constant
$C > 0$
such that
\[ v(c_1(L)+\epsilon\omega) \geq C \epsilon^{n-\nu(L)}. \]

\noindent
Let
$h = h_\infty e^{-2\psi}$
be a metric with
$\Theta_h(L) \geq 0$, 
and let
$\psi_\epsilon \downarrow \psi$
be an equisingular regularization of
$\psi$, 
such that
\[ \widetilde{h}_\epsilon := h_\infty e^{-2\psi_\epsilon} \]
satisfies 
$\Theta_{\widetilde{h}_\epsilon} \geq - \epsilon\omega$
in the sense of currents. The metrics considered in the follwing are given by
\[ \widehat{h}_\epsilon^{1+s} = 
   h_\infty \exp(-2(\delta({1+s})\phi_\epsilon + 
                    (1-\delta)({1+s})\psi_\epsilon)) \]
where
$\delta > 0$
is a sufficiently small number which will be fixed later. Note that
$\widehat{h}_\epsilon^{1+s}$
is really a metric on
$L$,
since
$h_\infty$
remains unchanged.

\noindent
$\widehat{h}_\epsilon^{1+s}$
is smooth outside an analytic subset
$Z_\epsilon \subset X$.
Its multiplier ideal is controlled by subadditivity 
(\cite[(14.2)]{Dem00}):
\[ \mathcal{J}(\widehat{h}_\epsilon^{1+s}) \subset 
   \mathcal{J}({h}_\epsilon^{\delta(1+s)}) \cdot
   \mathcal{J}(\widetilde{h}_\epsilon^{(1-\delta)(1+s)}) \subset 
   \mathcal{J}(\widetilde{h}_\epsilon^{(1-\delta)(1+s)}) =
   \mathcal{J}(h^{(1-\delta)(1+s)}) \]
because of the equisingularity. Locally, the H\"older inequality shows that
\[ \begin{array}{ll}
   \int |f|^2 e^{-2[\delta(1+s)\phi_\epsilon + 
                             (1-\delta)(1+s)\psi_\epsilon]}dV_\omega  \leq
    & \hspace{5cm}
   \end{array} \]
\[ \begin{array}{rr}
   \hspace{2cm} & (\int |f|^2 e^{-2(1+s)\phi_\epsilon} dV_\omega)^\delta 
        \cdot
        \int |f|^2  e^{-2(1+s)\psi_\epsilon} dV_\omega)^{1-\delta}, 
   \end{array} \]
hence 
\[ \mathcal{J}({h}_\epsilon^{1+s}) \cap \mathcal{J}(h^{1+s}) = 
   \mathcal{J}(h^{1+s}) \subset
   \mathcal{J}(\widehat{h}_\epsilon^{1+s}), \]
where the first equality comes from the properties of
$h_\epsilon$.

\noindent
Finally, we check how the eigenvalues of the curvature form are controlled: 
By construction,
\[ \begin{array}{rcl}
   \Theta_{\widehat{h}_\epsilon} + 2\epsilon\omega & = & 
   \delta(\Theta_{h_\epsilon}(L) + \epsilon\omega) + 
   (1-\delta)(\Theta_{\widetilde{h}_\epsilon}(L) + \epsilon\omega) 
   + \epsilon\omega \\
    & \geq & \delta(\Theta_{h_\epsilon}(L) + \epsilon\omega) + \epsilon\omega.
   \end{array}\]
At each point 
$x \in X \setminus Z_\epsilon$,
we may choose coordinate systems
$(z_j)_{1 \leq j \leq n}$
resp.
$(w_j)_{1 \leq j \leq n}$
which diagonalize simultaneously the hermitian forms
$\omega(x)$
and
$T_\epsilon+\epsilon\omega$
resp. 
$\Theta_{\widehat{h}_\epsilon} + 2\epsilon\omega$,
in such a way that
\[ \omega(x) = i \sum_{1 \leq j \leq n} dz_j \wedge d\overline{z}_j,\ 
   (T_\epsilon+\epsilon\omega)(x) = 
   i \sum_{1 \leq j \leq n} \lambda_j^{(\epsilon)} dz_j \wedge d\overline{z}_j
\]
resp.
\[ \omega(x) = i \sum_{1 \leq j \leq n} dw_j \wedge d\overline{w}_j,\ 
   (\Theta_{\widehat{h}_\epsilon} + 2\epsilon\omega)(x) = 
   i \sum_{1 \leq j \leq n} \widehat{\lambda}_j^{(\epsilon)} dw_j \wedge 
                                                           d\overline{w}_j. \]
Let
$\lambda_1^{(\epsilon)} \leq \ldots \leq \lambda_n^{(\epsilon)}$
and
$\widehat{\lambda}_1^{(\epsilon)} \leq \ldots \leq 
                                        \widehat{\lambda}_n^{(\epsilon)}$.
Changing from
$z_j$
to
$w_j$
by a unitary transformation, the
$\lambda_j^{(\epsilon)}$'s 
remain the same, and the inequality between the currents from above implies 
\[ \widehat{\lambda}_j^{(\epsilon)} \geq \delta \lambda_j^{(\epsilon)} + 
                                                                \epsilon, \]
by Weyl's monotonicity principle \cite[p.291]{Bha01}. On the other hand, the 
Monge-Amp\`ere inequality satisfied by 
$T_\epsilon$
tells us that
\[ \lambda_1^{(\epsilon)} \cdots \lambda_n^{(\epsilon)} \geq C \cdot 
   \epsilon^{n-\nu(L)} \]
almost everywhere on
$X$.

\section{Comparisons of the metrics}
\label{compmet-sec}

\noindent
Let
$Z_\epsilon$
be the analytic subset such that 
$\widehat{h}_\epsilon$
is smooth on
$X \setminus Z_\epsilon$.
By \cite[Prop.1.6]{Dem82}, for every 
$\epsilon > 0$
there is a sequence of complete K\"ahler metrics 
$(\omega_{\epsilon,t})$
on
$X \setminus Z_\epsilon$
converging from above against
$\omega_\epsilon = \omega$.

\noindent
Let
$\mathcal{D}^{n,q}_{c,\epsilon}$
be the space of all
$(n,q)$-\!
forms with values in
$L$
and coefficients in
$\mathcal{J}(\hat{h}_\epsilon^{1+s}) \otimes \mathcal{C}^\infty$
and compact support in
$X \setminus Z_\epsilon$.
Let
$L^{n,q}_{\epsilon,t}$
be the
$L^2$-\!
completion of
$\mathcal{D}^{n,q}_{c,\epsilon}$
with respect to the norm
\[ \parallel\! u \!\parallel^2_{\epsilon,t} := 
   \int_{X \setminus Z_\epsilon} 
   |u|_{\bigwedge^{n,q} \omega_{\epsilon,t} \otimes \hat{h}_\epsilon^{1+s}}^2
   dV_{\omega_{\epsilon,t}}, \]
including the case
$t = 0$,
where 
$\bigwedge^{n,q} \omega_{\epsilon,t} \otimes \hat{h}_\epsilon^{1+s}$
denotes the metric on 
$(n,q)$-\!
forms with values in
$L$
and coefficients in
$\mathcal{J}(\hat{h}_\epsilon^{1+s}) \otimes \mathcal{C}^\infty$
naturally induced by
$\omega_{\epsilon,t}$
and
$\hat{h}_\epsilon^{1+s}$.
The volume form
$dV_{\omega_{\epsilon,t}}$
equals
$\frac{\omega_{\epsilon,t}^n}{n!}$.

\noindent
The operator
$\overline{\partial}$
defines a linear, closed, densely defined operator 
\[ \overline{\partial}_{\epsilon,t}: 
   L^{n,q}_{\epsilon,t} \rightarrow L^{n,q+1}_{\epsilon,t}. \]
An element
$u \in L^{n,q}_{\epsilon,t}$
is in the domain
$D_{\overline{\partial}_{\epsilon,t}}$
if 
$\overline{\partial}(u)$,
defined in the sense of distribution theory, belongs to
$L^{n,q+1}_{\epsilon,t}$.
That
$\overline{\partial}_{\epsilon,t}$
is closed follows from the fact that differentiation is a continuous operation
in distribution theory, and the domain is dense since it contains
$\mathcal{D}^{n,q}_{c,\epsilon}$.

\noindent
Since
$\overline{\partial}_{\epsilon,t}$
is densely defined, there is an adjoint operator
$\overline{\partial}^\ast_{\epsilon,t}$,
and because 
$\overline{\partial}_{\epsilon,t}$
is closed,
\[ (\overline{\partial}^\ast_{\epsilon,t})^\ast = 
                             \overline{\partial}_{\epsilon,t})^\ast. \] 
(cf. \cite[p.29]{Nag67}). Let
$D_{\overline{\partial}_{\epsilon,t}^\ast}$
denote the domain of the operator
$\overline{\partial}_{\epsilon,t}$
in
$L^{n,q}_{\epsilon,t}$.

\noindent
Let
$\Lambda$
be the adjoint of the operator
$L$
which multiplicates with
$\omega$,
that is
\[ L\alpha = \omega \wedge \alpha,\ \ \ 
   \langle \Lambda\alpha|\beta \rangle_{\epsilon,t} = 
   \langle \alpha|\ \omega \wedge \beta \rangle_{\epsilon,t} \]
for all forms
$\alpha, \beta \in L^{n,q}_{\epsilon,t}$.
(The scalar product above is taken in every point
$z \in X \setminus Z_\epsilon$.)

\noindent
If
$\theta$
is a real
$(1,1)$-\!
form we define for all 
$q = 1, \ldots, n$
a sesquilinear form
$\theta_q$
on the fibers of
$\Omega^{n,q}_{X \setminus Z_\epsilon} \otimes L$
by setting in every point
$z \in X \setminus Z_\epsilon$
\[ \theta_q(\alpha,\beta) = 
   \langle \theta\Lambda\alpha|\beta \rangle_{\epsilon,t} \]
for all
$\alpha, \beta \in \Omega^{n,q}_{X \setminus Z_\epsilon,z} \otimes L_z$.
If
$\theta = \Theta_{\epsilon,t}$
is the curvature form of the metric 
$\hat{h}_\epsilon^{1+s}$
on
$L$,
the term
$\langle \theta\Lambda\alpha|\beta \rangle_{\epsilon,t}$
occurs in the Bochner-Kodaira inequality:
\[ \parallel\! \overline{\partial}_{\epsilon,t} u \!\parallel^2_{\epsilon,t} + 
   \parallel\! \overline{\partial}_{\epsilon,t}^\ast u 
                                                  \!\parallel^2_{\epsilon,t} 
   \geq \int_X  \langle \Theta_{\epsilon,t}\Lambda u|u \rangle_{\epsilon,t}
   dV_{\omega_{\epsilon,t}}. \]
On
$\mathcal{D}^{n,q}_{c,\epsilon}$,
this inequality is valid by the usual computations 
(\cite[(4.7)]{Dem00}). H\"ormander (\cite[Lem. 5.2.1]{Hoer65}) showed that for
the complete metric
$\omega_{\epsilon,t}$
($t > 0$)
the forms in
$\mathcal{D}^{n,q}_{c,\epsilon}$
are dense in
$D_{\overline{\partial}_{\epsilon,t}} \cap 
 D_{\overline{\partial}_{\epsilon,t}^\ast}$
w.r.t. the graph norm
\[ u \mapsto \parallel\! u \!\parallel^2_{\epsilon,t} + 
   \parallel\! \overline{\partial}_{\epsilon,t} u \!\parallel^2_{\epsilon,t} + 
   \parallel\! \overline{\partial}^\ast_ {\epsilon,t}u 
                                                  \!\parallel^2_{\epsilon,t}.\]
Hence we have the Bochner inequality for all
$u \in D_{\overline{\partial}_{\epsilon,t}} \cap 
       D_{\overline{\partial}_{\epsilon,t}^\ast}$.

\noindent
To really apply the Bochner technique we still need some comparative 
inequalities between the different metrics
$\omega_{\epsilon,t}$:

\begin{lem} \label{easyineq-lem}
$\parallel\! u \!\parallel_{\epsilon,t^\prime} \leq 
 \parallel\! u \!\parallel_{\epsilon,t}$
for all
$u \in L^{n,q}_{\epsilon,t}$
and all
$0 \leq t \leq t^\prime $.
\end{lem}
\begin{proof}
This is just Lemma 3.3 in \cite{Dem82}. 
\end{proof}

\noindent
Consequently, we have a linear continuous operator
$f_q: L^{n,q}_{\epsilon,t} \rightarrow L^{n,q}_{\epsilon,t^\prime}$
with norm
$\parallel\! f_q \!\parallel \leq 1$.
\begin{lem} \label{commdia-lem}
The diagram 
\[ \xymatrix{\ar @{} [dr] |{\#} 
   L^{n,q}_{\epsilon,t} \ar[r]^{f_q} 
                            \ar[d]_{\overline{\partial}_{\epsilon,t}} & 
   L^{n,q}_{\epsilon,t^\prime} 
                     \ar[d]^{\overline{\partial}_{\epsilon,t^\prime}} \\
   L^{n,q+1}_{\epsilon,t} \ar[r]_{f_{q+1}} & L^{n,q+1}_{\epsilon,t^\prime}
   } \]
is commutative.
\end{lem}
\begin{proof}
Let
$v_{\epsilon,t}$
be any element of
$L^{n,q}_{\epsilon,t}$
such that
$\overline{\partial}_{\epsilon,t}$
is defined. Since
$\mathcal{D}^{n,q}_{c,\epsilon}$
is dense in
$D_{\overline{\partial}_{\epsilon,t}}$
there is a sequence of smooth forms
$(v_{\epsilon,t}^{(n)})_{n \in \mathbb{N}}$
in
$\mathcal{D}^{n,q}_{c,\epsilon}$
such that
\[ v_{\epsilon,t}^{(n)} \rightarrow v_{\epsilon,t},\ \ 
   \overline{\partial} v_{\epsilon,t}^{(n)} \rightarrow 
   \overline{\partial}_{\epsilon,t} v_{\epsilon,t} \]
strongly in the 
$(\epsilon,t)$-\!
norm. Now, the derivation of smooth forms in 
$\mathcal{D}^{n,q}_{c,\epsilon}$
w.r.t.
$\overline{\partial}_{\epsilon,t}$
and
$\overline{\partial}_{\epsilon,t^\prime}$ 
do not differ. So the two limits above exist and are the same for the
$(\epsilon,t^\prime)$-\!
norm because of lemma~\ref{easyineq-lem}.
\end{proof}

\noindent
Again, let
$\theta$
be any real
$(1,1)$-\!
form. If
$\alpha \in L^{n,q}_{\epsilon,t}$
we define 
$|\alpha|_\theta$
in every point
$z \in X \setminus Z_\epsilon$
as the smallest number 
$\geq 0$
(perhaps infinte) such that
\[ |\langle \alpha|\beta \rangle_{\epsilon,t}|^2 \leq
   |\alpha|_\theta^2  \langle \theta\Lambda\alpha|\beta \rangle_{\epsilon,t}\]
for all
$\beta \in L^{n,q}_{\epsilon,t}$.
\begin{lem} \label{normineq-lem}
The 
$(n,n)$-\!
form
$|\alpha|_\theta^2dV_{\omega_{\epsilon,t}}$
decreases if 
$t$
increases.
\end{lem}
\begin{proof}
This is just Lemma 3.2 in \cite{Dem82}.
\end{proof} 

\begin{lem} \label{diffineq-lem}
For all
$\beta \in L^{n,q}_{\epsilon,t}$,
we have
\[ \int_X |\beta|_\theta^2 dV_{\omega_{\epsilon,t}} \leq 
   \int_X \frac{1}{\lambda_1 + \ldots + \lambda_q} 
          |\beta|^2_{\omega_{\epsilon,t}} dV_{\omega_{\epsilon,t}}, \]
where
$\lambda_1, \ldots, \lambda_q$
are the 
$q$
smallest eigenvalues of
$\theta$
with respect to
$\omega_{\epsilon,t}$.
\end{lem}
\begin{proof}
There is an orthonormal base 
$dz_1, \ldots, dz_n$
of
$\Omega^{n,q}_{X \setminus Z_\epsilon,z}$
such that we can write
\[ \omega_{\epsilon,t} = \frac{i}{2} \sum_{j=1}^n dz_j \wedge 
                                                        d\overline{z}_j, \]
\[ \theta = \frac{i}{2} \sum_{j=1}^n \lambda_j dz_j \wedge 
                              d\overline{z}_j,\ \ \lambda_j \in \mathbb{R}. \]
$\beta \in L^{n,q}_{\epsilon,t}$
may be written as
\[ \beta = \sum_{|J|=q} \beta_J dz_1 \wedge \ldots \wedge dz_n \wedge 
           d\overline{z}_J \otimes e, \]
where
$e$
is any section in
$L$
which makes the
$dz_1 \wedge \ldots \wedge dz_n \wedge 
           d\overline{z}_J \otimes e $
orthonormal in
$L^{n,q}_{\epsilon,t}$. 
We verify
\[ \Lambda \beta = 2 \sum_{|J|=q-1} \sum_{1 \leq j \leq n} (-1)^{n-j}
   \beta_{jJ} dz_1 \wedge \ldots \wedge \widehat{dz_j} \wedge \ldots 
                   \wedge dz_n \wedge d\overline{z}_J \otimes e \]
($\widehat{dz_j}$ 
meaning that we omit
$dz_j$), 
and
\[ \theta_q(\beta,\beta) = 
   \langle \theta\Lambda\beta|\beta \rangle_{\epsilon,t} = 
   2^{n+q} \sum_{|J|=q-1} \sum_{1 \leq j \leq n} \lambda_j \beta_{jJ}
                                                 \overline{\beta_{jJ}} = 
   2^{n+q} \sum_{|J|=q} (\sum_{j \in J} \lambda_j) |\beta_J|^2. \]
Now,
\[ \begin{array}{rcl} 
   |\beta|_\theta^2 & = & 
   \sup_u \frac{|\langle \beta|\ u \rangle_{\epsilon,t}|^2}
               {\langle \theta\Lambda u|\ u \rangle_{\epsilon,t}} \leq 
   \sup_u \frac{|\beta|_{\epsilon,t}^2 \cdot |u|_{\epsilon,t}^2}
               {\langle \theta\Lambda u|\ u \rangle_{\epsilon,t}} = 
   |\beta|_{\epsilon,t}^2 \cdot \sup_u \frac{|u|_{\epsilon,t}^2}
                     {\langle \theta\Lambda u|\ u \rangle_{\epsilon,t}} = \\
    &   & \\
    & = & |\beta|_{\epsilon,t}^2 \cdot \sup_u \frac{\sum_{|J|=q} |u_J|^2}
                   {\sum_{|J|=q} (\sum_{j \in J} \lambda_j) |u_J|^2} \leq
          |\beta|_{\epsilon,t}^2 \frac{1}{\lambda_1 + \ldots + \lambda_q},
   \end{array} \]
where 
$u = \sum_{|J|=q} u_J dz_1 \wedge \ldots \wedge dz_n \wedge 
           d\overline{z}_J \otimes e$
as above. The lemma follows.
\end{proof}

\section{The Bochner technique}
\label{Bochner-sec}

\noindent
Let
$K^{n,q}_{\epsilon,t}$
be the 
$L^2$-\!
completion of
$\mathcal{D}^{n,q}_{c,\epsilon}$
w.r.t. the metric
\[ \parallel\! u \!\parallel^2_{K^{n,q}_{\epsilon,t}} := 
   \int_{X \setminus Z_\epsilon} 
   (|u|_{\bigwedge^{n,q} \omega_{\epsilon,t} \otimes 
    \hat{h}_\epsilon^{1+s}}^2 + 
    |\overline{\partial} u|_{\bigwedge^{n,q} \omega_{\epsilon,t} \otimes 
    \hat{h}_\epsilon^{1+s}}^2)
   dV_{\omega_{\epsilon,t}}, \]
in the space of all forms with
$L^2_{\mathrm{loc}}$
coefficients, and let
$\mathcal{K}^{n,q}_{\epsilon,t}$
be the corresponding sheaf of germs of locally
$L^2$
sections on
$X$
(the local
$L^2$
condition should hold on
$X$
and not only on
$X \setminus Z_\epsilon$). 
\begin{lem} \label{exact-lem} 
For all
$\epsilon > 0$,
the
$L^2$
Dolbeault complex
$(\mathcal{K}^{n,q}_{\epsilon,0}, \overline{\partial}_{\epsilon,0})$
is a fine resolution of the sheaf
$K_X \otimes L \otimes \mathcal{J}(\hat{h}_\epsilon^{1+s})$.
\end{lem}
\begin{proof}
$X \setminus Z_\epsilon$
can be covered by open subsets
$U \subset X$
which are Stein, lie relatively compact in another Stein open subset
$U^\prime \subset X$, 
and on which 
$L$
is trivial. On these
$U$'s 
we can show the
$\overline{\partial}$-
Poincar\'e lemma.

\noindent
First, as Stein sets,
$U \subset\subset U^\prime$
may be embedded as analytic subsets into some
$\mathbb{C}^N$. 
Hence we can find a smooth plurisubharmonic function
$\psi$
on
$U^\prime$
such that
$i\partial\overline{\partial}\psi \geq 2\lambda\omega$
for some constant
$\lambda > 0$
on
$U$
($\omega_{\epsilon,0} = \omega$
is smooth on
$U$). 
Furthermore,
$|\psi|$
is bounded on
$U$
by some constant
$M > 0$. 
Subtracting 
$M$
we still have a plurisubharmonic function, which we also call
$\psi$,
satisfying
\[ -2M \leq \psi \leq 0\ \ \ \mathrm{and}\ \ \ 
   i\partial\overline{\partial}\psi \geq 2\lambda\omega. \]
This implies that the metrics
$\widehat{h}_\epsilon^{1+s}$
and
$\widehat{h}_\epsilon^{1+s} e^{-2\psi}$
are comparable on 
$U$, 
and
$\Theta_{\widehat{h}_\epsilon^{1+s} e^{-2\psi}} \geq \lambda\omega$,
for 
$\lambda$
sufficiently big.

\noindent
Since
$\widehat{h}_\epsilon^{1+s}$
and
$\widehat{h}_\epsilon^{1+s} e^{-2\psi}$
are comparable, we can interchange them in the metric for
$\mathcal{K}^{n,q}_{\epsilon,0}(U)$.
So, given
$g \in \mathcal{K}^{n,q}_{\epsilon,0}(U)$
with
$\overline{\partial}_{\epsilon,0}(g) = 0$,
we know that 
$\int_X |g|^2_{\epsilon,0} dV_\omega < \infty$,
hence also
$\int_X |g|^2_{\epsilon,\psi} dV_\omega < \infty$
with the new metric, and
\[ \int_X |g|^2_{\Theta_{\widehat{h}_\epsilon^{1+s} e^{-2\psi}}} dV_\omega \leq
   \int_X \frac{1}{q\lambda}|g|^2_{\epsilon,\psi} dV_\omega < \infty. \]
Therefore, we can apply theorem 4.1 of \cite{Dem82}: There exists a 
$(n,q-1)$-\!
form
$f$
with
$L^2_{\mathrm{loc}}$
coefficients in
$U$
such that 
\[ \overline{\partial}_{\epsilon,0}(f) = g \]
and
\[ \int_X |f|^2_{\epsilon,0} dV_\omega \leq 
   \int_X |f|^2_{\epsilon,\psi} dV_\omega \leq 
   \int_X \frac{1}{q\lambda}|g|^2_{\epsilon,\psi} dV_\omega \leq
   e^{2M} \int_X |g|^2_{\epsilon,0} dV_\omega. \]

\noindent
Finally, the 
$L^2$
condition forces sections holomorphic on
$X \setminus Z_\epsilon$
to extend holomorphically across
$Z_\epsilon$
(\cite[Lem.6.9]{Dem82}). The 
$L^2$
condition implies that the coefficients lie in
$\mathcal{J}(\hat{h}_\epsilon^{1+s})$.
Consequently, the complex is a resolution, and it is 
fine because of the existence of partition of unities.
\end{proof} 

\noindent
Now, let us take a cohomology class 
$\{\beta\} \in H^q(X, \mathcal{O}(K_X + L) \otimes 
    \mathcal{J}(h_{\min}^{1+s^\prime}))$.
If
$\mathcal{U}$
is a covering of
$X$
with Stein open subsets
$U_\alpha$,
the class
$\{\beta\}$
may be represented by a \v{C}ech cocycle
\[ (\beta_{\alpha_0 \cdots \alpha_q})_{\alpha_0 \cdots \alpha_q} \in
   C^q(\mathcal{U}, \mathcal{O}(K_X + L) \otimes 
                                  \mathcal{J}(h_{\min}^{1+s^\prime})) \subset
   C^q(\mathcal{U}, \mathcal{O}(K_X + L) \otimes 
                               \mathcal{J}(\hat{h}_{\epsilon}^{1+s})). \]
Let
$(\psi_\alpha)$
be a
$\mathcal{C}^\infty$
partition of unity subordinate to
$\mathcal{U}$.
Taking the usual De Rham-Weil isomorphisms between \v{C}ech and Dolbeault 
cohomology, we obtain a 
closed
$(n,q)$-\!
form in
$K^{n,q}_{\epsilon,0}$
of the form
\[ \beta = \sum_{\alpha_0, \ldots, \alpha_q} \beta_{\alpha_0 \cdots \alpha_q}
           \overline{\partial} \psi_{\alpha_0} \wedge \ldots \wedge
           \overline{\partial} \psi_{\alpha_q}. \]
In particular, this form has coefficients in
$\mathcal{J}(\hat{h}_{\epsilon}^{1+s}) \otimes \mathcal{C}^\infty$.
We want to show that
$\beta$
is a boundary in
$K^{n,q}_{\epsilon,0}$
for some
$\epsilon > 0$,
hence
$\{\beta\} = 0 \in
 H^q(X, \mathcal{O}(K_X + L) \otimes \mathcal{J}(\hat{h}_{\epsilon}^{1+s}))$.
This implies theorem~\ref{KV0-thm} because of the inclusion
$\mathcal{J}(\hat{h}_{\epsilon}^{1+s}) \subset
 \mathcal{J}(h_{\min}^{(1-\delta)(1+s)})$.

\noindent
The reasoning starts as follows: 
$\beta$
is also an element of
$K^{n,q}_{\epsilon,t}$
for any
$t \geq 0$,
because of lemma~\ref{easyineq-lem}.
Every 
$L^2$
form
$u \in  D_{\overline{\partial}_{\epsilon,t}^\ast} \subset 
                                                       K^{n,q}_{\epsilon,t}$ 
may be written as
$u = u_1 + u_2$
with
\[ u_1 \in \mathrm{ker\ } \overline{\partial}_{\epsilon,t}\ \ \mathrm{and}\ \ 
   u_1 \in (\mathrm{ker\ } \overline{\partial}_{\epsilon,t})^\perp = 
           \overline{\mathrm{im\ } \overline{\partial}_{\epsilon,t}^\ast} 
   \subset \mathrm{ker\ } \overline{\partial}_{\epsilon,t}^\ast, \]
since
$\overline{\partial}_{\epsilon,t}$
is a closed operator, hence
$\mathrm{ker\ } \overline{\partial}_{\epsilon,t}$
is closed. Using
$\beta \in \mathrm{ker\ } \overline{\partial}$
and the two inequalities in lemma~\ref{normineq-lem} and \ref{diffineq-lem}, we
get 
($\Theta_{\epsilon,t}$
denotes the curvature form of
$\hat{h}_{\epsilon}^{1+s}$
on
$X \setminus Z_\epsilon$,
plus
$2\epsilon\omega_{\epsilon,t}$)
\[ \begin{array}{rcl}
   |\ll \beta,u \gg_{\epsilon,t}|^2 & = &
   |\ll \beta,u_1 \gg_{\epsilon,t}|^2 = 
   |\int_{X \setminus Z_\epsilon} \langle \beta,u_1 \rangle_{\epsilon,t} 
                                                dV_{\omega_{\epsilon,t}}|^2 
      \leq  \\
    &      & \\
    & \leq & (\int_{X \setminus Z_\epsilon} 
             |\langle \beta,u_1 \rangle_{\epsilon,t}| 
             dV_{\omega_{\epsilon,t}})^2 \leq \\
    &      & \\
    & \leq & (\int_{X \setminus Z_\epsilon} \beta_{\Theta_{\epsilon,t}} \cdot
             \sqrt{\langle \Theta_{\epsilon,t}\Lambda u_1|u_1 
                                                   \rangle_{\epsilon,t}} 
             dV_{\omega_{\epsilon,t}})^2 \leq \\
    &      & \\
    & \leq & \int_{X \setminus Z_\epsilon} \beta_{\Theta_{\epsilon,t}}^2
             dV_{\omega_{\epsilon,t}} \cdot \int_{X \setminus Z_\epsilon} 
             \langle \Theta_{\epsilon,t}\Lambda u_1|u_1 \rangle_{\epsilon,t} 
             dV_{\omega_{\epsilon,t}} \leq \\
    &      & \\
    & \leq & \int_{X \setminus Z_\epsilon} \beta_{\Theta_{\epsilon,0}}^2
             dV_{\omega_{\epsilon,0}} \cdot \int_{X \setminus Z_\epsilon}
             \langle \Theta_{\epsilon,t}\Lambda u_1|u_1 \rangle_{\epsilon,t} 
             dV_{\omega_{\epsilon,t}} \leq \\
    &      & \\
    & \leq & \int_{X \setminus Z_\epsilon} 
             \frac{1}{\hat{\lambda}^{(\epsilon,0)}_1 + \cdots + 
                   \hat{\lambda}^{(\epsilon,0)}_q} |\beta|_{\epsilon,0}^2
                                                 dV_{\omega_{\epsilon,0}} \cdot
             \int_{X \setminus Z_\epsilon}
             \langle \Theta_{\epsilon,t}\Lambda u_1|u_1 \rangle_{\epsilon,t}
             dV_{\omega_{\epsilon,t}}. 
   \end{array} \]

\noindent
$u_1$
is an element of
$D_{\overline{\partial}_{\epsilon,t}} \cap 
       D_{\overline{\partial}_{\epsilon,t}^\ast}$,
since
$u_1 \in \mathrm{ker\ } \overline{\partial}_{\epsilon,t}$,
$u_2 \in \mathrm{ker\ } \overline{\partial}_{\epsilon,t}^\ast$
and 
$u_1 = u - u_2$.
Consequently, we can apply the Bochner inequality on
$u_1$.
As
$\overline{\partial} u_1 = 0$
we get that the second integral on the right hand side is bounded above by
\[ \parallel\! \overline{\partial}^\ast_{\epsilon,t} u_1 
                                             \!\parallel^2_{\epsilon,t} +\  
   2q\epsilon\!\parallel\! u_1 \!\parallel^2_{\epsilon,t}\ \leq\ 
   \parallel\! \overline{\partial}^\ast_{\epsilon,t} u 
                                                 \!\parallel^2_{\epsilon,t} +\ 
   2q\epsilon\!\parallel\! u \!\parallel^2_{\epsilon,t}, \]
and finally
\[ |\langle \beta,u \rangle_{\epsilon,t}|^2 \leq 
   \int_X \frac{1}{\hat{\lambda}^{(\epsilon,t)}_1 + \cdots + 
                   \hat{\lambda}^{(\epsilon,t)}_q} |\beta|_{\epsilon,t}^2
                                                 dV_{\omega_{\epsilon,t}} 
   (\parallel\! \overline{\partial}^\ast_{\epsilon,t} u 
                                                 \!\parallel^2_{\epsilon,t} +\ 
   2q\epsilon\!\parallel\! u \!\parallel^2_{\epsilon,t}), \]
where the term
$2q\epsilon\!\parallel\! u \!\parallel^2_{\epsilon,t})$
comes in because
$\Theta_{\epsilon,t}$
differs from the curvature form of
$\hat{h}_{\epsilon}^{1+s}$
by
$2\epsilon\omega_{\epsilon,t}$.

\noindent
Using the uniform bound
$C_\epsilon = \int_X \frac{1}{\hat{\lambda}^{(\epsilon,0)}_1 + \cdots + 
                   \hat{\lambda}^{(\epsilon,0)}_q} |\beta|_{\epsilon,0}^2
                                                 dV_{\omega_{\epsilon,0}}$
we apply the Hahn-Banach theorem as in~\cite{Dem82}:
$\ll\! \beta,u \gg_{\epsilon,t}$
defines a linear form on the range of the densely defined operator
\[ T: L^{n,q}_{\epsilon,t} \rightarrow L^{n,q-1}_{\epsilon,t} \oplus 
                                       L^{n,q}_{\epsilon,t},\ \ 
   u \mapsto \overline{\partial}^\ast_{\epsilon,t} u + 2q\epsilon u \]
(with domain
$D_T = D_{\overline{\partial}^\ast_{\epsilon,t}}$). 
Hence there exists an
$f_{\epsilon,t} = v_{\epsilon,t} \oplus \frac{1}{2q\epsilon}w_{\epsilon,t}
   \in L^{n,q-1}_{\epsilon,t} \oplus L^{n,q}_{\epsilon,t}$
such that
\[ \ll \beta,u \gg_{\epsilon,t} = 
   \ll f_{\epsilon,t}, \overline{\partial}^\ast_{\epsilon,t} u + 
                          2q\epsilon u \gg_{\epsilon,t}. \]
Consequently,
$\beta = T^\ast f_{\epsilon,t} = 
         \overline{\partial}_{\epsilon,t} v_{\epsilon,t} + w_{\epsilon,t}$
with
\[ \parallel\! v_{\epsilon,t} \!\parallel^2_{\epsilon,t} +\ 
   \frac{1}{2q\epsilon}\parallel\! w_{\epsilon,t} \!\parallel^2_{\epsilon,t}
   \leq C_\epsilon. \]
Furthermore,
$ \overline{\partial}_{\epsilon,t} w_{\epsilon,t} = \overline{\partial} \beta 
  = 0$, 
and
$v_{\epsilon,t}, w_{\epsilon,t}$
are both contained in
$K^{n,q}_{\epsilon,t}$.

\noindent 
The estimates of section~\ref{compmet-sec} tell us that the metrics of the
$v_{\epsilon,t}$
and
$w_{\epsilon,t}$
in the 
$L^2$ 
space
$L^{n,q-1}_{\epsilon,t_0}$
resp.
$L^{n,q}_{\epsilon,t_0}$
are uniformly bounded for all
$t_0 \geq t > 0$.
Since
$\parallel\! \beta \!\parallel_{\epsilon,t_0} \geq
 |\parallel\! \overline{\partial}_{\epsilon,t} v_{\epsilon,t} 
                                               \!\parallel_{\epsilon,t_0} -
  \parallel\! w_{\epsilon,t} \!\parallel_{\epsilon,t_0}|$,
the same is true for
$\overline{\partial}_{\epsilon,t} v_{\epsilon,t}$.
Consequently, the sequences of these elements converge weakly as 
$t \rightarrow 0$, 
and we have three limits in the respective spaces:
\[ v_{\epsilon,t} \rightharpoonup v_\epsilon \in L^{n,q-1}_{\epsilon,t_0},\ \ 
   \overline{\partial}_{\epsilon,t} v_{\epsilon,t} \rightharpoonup 
   v_\epsilon^\prime \in L^{n,q}_{\epsilon,t_0},\ \ 
   w_{\epsilon,t} \rightharpoonup w_\epsilon \in L^{n,q}_{\epsilon,t_0}. \]
Note that these weak limits are the same for every choice of
$t_0 > 0$: 
The spaces
$L^{n,q}_{\epsilon,t_0}$
all contain the dense subset
$\mathcal{D}^{n,q}_{c,\epsilon}$,
hence weak convergence is transmitted through the continuous maps between 
them. 

\vspace{0.2cm}

\noindent
\textit{Claim}: 
$\overline{\partial}_{\epsilon,t_0} v_\epsilon = v_\epsilon^\prime$.
\begin{proof}
On the one hand, we have
\[ \ll \overline{\partial}_{\epsilon,t} v_{\epsilon,t}, 
                                       u \gg_{\epsilon,t_0} \rightarrow
   \ll v_\epsilon^\prime, u \gg_{\epsilon,t_0} \]
for all 
$u \in D_{\overline{\partial}^\ast_{\epsilon,t_0}}$ 
because of the weak convergence.
On the other hand, the commutativity of the diagram in lemma~\ref{commdia-lem}
and again the weak convergence show that
\[  \begin{array}{ll}
    \ll \overline{\partial}_{\epsilon,t} v_{\epsilon,t}, 
                                       u \gg_{\epsilon,t_0} =
    \ll \overline{\partial}_{\epsilon,t_0} v_{\epsilon,t}, 
                                       u \gg_{\epsilon,t_0} =
    \ll v_{\epsilon,t}, \overline{\partial}^\ast_{\epsilon,t_0} 
                                          u \gg_{\epsilon,t_0} \rightarrow &
    \hspace{4cm}
    \end{array} \]
\[ \begin{array}{rr}
    \hspace{5cm} &
    \ll v_\epsilon, \overline{\partial}^\ast_{\epsilon,t_0} 
                                          u \gg_{\epsilon,t_0} = 
    \ll \overline{\partial}_{\epsilon,t_0} v_\epsilon, 
                                                 u \gg_{\epsilon,t_0}. 
    \end{array} \]
Since
$D_{\overline{\partial}^\ast_{\epsilon,t_0}}$
is dense in
$L^{n,q}_{\epsilon,t_0}$
the claim follows. 
\end{proof}

\noindent
By standard properties of weak convergence,
\[ \parallel\! v_\epsilon \!\parallel_{\epsilon,t_0} \leq
   \liminf_{t \rightarrow 0} \parallel\! v_{\epsilon,t} 
                                               \!\parallel_{\epsilon,t_0}, \]
and similarly for
$\parallel\! \overline{\partial}_{\epsilon,t_0} v_\epsilon 
                                               \!\parallel_{\epsilon,t_0}$
and
$\parallel\! w_\epsilon \!\parallel_{\epsilon,t_0}$. 
Consequently, these three norms are uniformly bounded by 
$C_\epsilon$
for
$t_0 > 0$.

\noindent
Now, we restrict the integral defining the 
$(\epsilon,t_0)$-\!
norm to compact subsets
$K \subset X \setminus Z_\epsilon$.
Of course, we get
$\parallel\! v_\epsilon \!\parallel_{\epsilon,t_0,K} \leq 
 \parallel\! v_\epsilon \!\parallel_{\epsilon,t_0}$,
hence the new
$(\epsilon,t_0,K)$-\!
norms of 
$v_\epsilon$
are still uniformly bounded by
$C_\epsilon$
in
$t_0 > 0$.
Furthermore, as
$\omega_t \downarrow \omega$,
we see that 
$\parallel\! v_\epsilon \!\parallel_{\epsilon,t_0,K} \rightarrow
 \parallel\! v_\epsilon \!\parallel_{\epsilon,0,K}$,
and monotone convergence tells us that
$\parallel\! v_\epsilon \!\parallel_{\epsilon,0}$
exists and is
$\leq C_\epsilon$.
The same is true for
$\parallel\! \overline{\partial}_{\epsilon,0} v_\epsilon 
                                               \!\parallel_{\epsilon,0}$
and
$\parallel\! w_\epsilon \!\parallel_{\epsilon,0}$.

\noindent
As
$\beta = \overline{\partial}_{\epsilon,t} v_{\epsilon,t} + w_{\epsilon,t}$
for all
$t$,
$\beta = \overline{\partial}_{\epsilon,0} v_\epsilon + w_\epsilon$
remains true. Furthermore, 
$\overline{\partial}_{\epsilon,0} w_{\epsilon,t} = 
 \overline{\partial}_{\epsilon,t} w_{\epsilon,t} = 0$.
Hence,
$v_\epsilon$
and
$w_\epsilon$
belong to
$K^{n,q-1}_{\epsilon,0}$
resp. 
$K^{n,q}_{\epsilon,0}$.

\noindent
For the last step we note that the almost plurisubharmonic weights
$\psi_\epsilon$
defining
$\tilde{h}_\epsilon$
form a decreasing sequence, and consequently,
\[ \parallel\! u \!\parallel_{\epsilon,0} \leq 
   \parallel\! u \!\parallel_{\epsilon^\prime,0}\ \ \ \ 
   \forall u \in K^{n,q}_{\epsilon^\prime,0}, \]
if
$\epsilon^\prime < \epsilon$.
This implies 
$\parallel\! w_\epsilon \!\parallel_{\epsilon_0,0} \leq
 \parallel\! w_\epsilon \!\parallel_{\epsilon,0}$
for some fixed
$\epsilon_0 > 0$
and
$\epsilon < \epsilon_0$.
Since
\[ \parallel\! w_\epsilon \!\parallel_{\epsilon,0} \leq C_\epsilon = 
   2q\epsilon \cdot C_\epsilon, \]
we conclude with the estimates in section~\ref{uniesti-sec} that
$\parallel\! w_\epsilon \!\parallel_{\epsilon_0,0} \rightarrow 0$
for
$\epsilon \rightarrow 0$.
But the norm of 
$w_\epsilon$
measures the distance of
$\beta$
from the closure of the subspace of boundaries in
$K^{n,q}_{\epsilon_0,0}$.
So it only remains to show
\begin{lem}
The subspace 
$B^{n,q}_{\epsilon_0,0} \subset K^{n,q}_{\epsilon_0,0}$
of boundaries in the Dolbeault complex 
$(\mathcal{K}^{n,q}_{\epsilon_0,0}, \overline{\partial})$
is closed.
\end{lem}
\begin{proof}
Let
$Z^{n,q}_{\epsilon_0,0} \subset K^{n,q}_{\epsilon_0,0}$
be the space of cocycles with respect to
$\overline{\partial}$,
and let
$\mathcal{Z}^{n,q}_{\epsilon_0,0} \subset \mathcal{K}^{n,q}_{\epsilon_0,0}$
be the corresponding sheaf. Let
$\mathcal{U}$
be a covering of 
$X$
with Stein open subsets as in the proof of exactness of the Dolbeault complex
in lemma~\ref{exact-lem}. By the usual DeRham-Weil isomorphism, 
\[ H^q(K^{n,\bullet}_{\epsilon_0,0}) = \frac{Z^{n,q}_{\epsilon_0,0}}
                             {\overline{\partial} K^{n,q-1}_{\epsilon_0,0}} =
   \frac{Z^0(\mathcal{U},\mathcal{Z}^{n,q}_{\epsilon_0,0})}
  {\overline{\partial} Z^0(\mathcal{U},\mathcal{K}^{n,q-1}_{\epsilon_0,0})}. \]
So we have to prove that
$\overline{\partial} Z^0(\mathcal{U},\mathcal{K}^{n,q}_{\epsilon_0,0})$
is closed in
$Z^0(\mathcal{U},\mathcal{Z}^{n,q}_{\epsilon_0,0})$
with respect to the 
$L^2$
norms on every set
$U$
in
$\mathcal{U}$.

\noindent
Note first that 
\[ \overline{\partial}: C^0(\mathcal{U},\mathcal{K}^{n,q-1}_{\epsilon_0,0})
   \rightarrow C^0(\mathcal{U},\mathcal{Z}^{n,q}_{\epsilon_0,0}) \]
is continuous, by definition of the norms, and surjective, by exactness. Hence 
$\overline{\partial}$
is an open map, by Banach's open mapping theorem. Its kernel is
$C^0(\mathcal{U},\mathcal{Z}^{n,q-1}_{\epsilon_0,0})$.

\noindent
Next, 
$C^0(\mathcal{U},\mathcal{Z}^{n,q-1}_{\epsilon_0,0})$
and
$Z^0(\mathcal{U},\mathcal{K}^{n,q-1}_{\epsilon_0,0})$
are closed in
$C^0(\mathcal{U},\mathcal{K}^{n,q-1}_{\epsilon_0,0})$,
since
$\overline{\partial}$
is a closed operator, and equality is conserved when going to the limit. 
Consequently, the sum of these two spaces is closed, too, and its complement
is open. But then
\[ \overline{\partial}(Z^0(\mathcal{U},\mathcal{K}^{n,q-1}_{\epsilon_0,0}) +
                       C^0(\mathcal{U},\mathcal{Z}^{n,q-1}_{\epsilon_0,0})) =
   \overline{\partial}(Z^0(\mathcal{U},\mathcal{K}^{n,q-1}_{\epsilon_0,0})) \]
is closed in
$C^0(\mathcal{U},\mathcal{Z}^{n,q}_{\epsilon_0,0})$
and also in
$Z^0(\mathcal{U},\mathcal{Z}^{n,q}_{\epsilon_0,0})$.
\end{proof}

\section{The uniform estimate}
\label{uniesti-sec}

\noindent
The aim of this section is to prove
\begin{lem}
Let
$0 < s < s^\prime$.
For every smooth
$(n,q)$-\!
form
$\beta$
with values in
$L$
and coefficients in
$\mathcal{J}(h^{1+s^\prime}_{\min}) \otimes \mathcal{C}^\infty$,
\[ \int_X \frac{q\epsilon}{\hat{\lambda}^{(\epsilon,0)}_1 + \cdots + 
                   \hat{\lambda}^{(\epsilon,0)}_q} |\beta|_{\epsilon,0}^2
                                               dV_{\omega_{\epsilon,0}},\]
tends to
$0$
for
$\epsilon \rightarrow 0$.
\end{lem}

\noindent
Before starting with the proof, set
$\hat{\lambda}_j := \hat{\lambda}^{(\epsilon,0)}_j$
and
$\lambda_j := \lambda^{(\epsilon,0)}_j$.

\noindent
Now, by construction we know that
$\hat{\lambda}_j \geq \delta \lambda_j + \epsilon$,
and
\[ \lambda_q^q\lambda_{q+1} \ldots \lambda_n \geq \lambda_1 \ldots \lambda_n
   \geq C \epsilon^{n-\nu}, \]
hence
\[ \frac{1}{\lambda_1 + \ldots + \lambda_q} \leq \frac{1}{\lambda_q} \leq
   C^{-1/q} \epsilon^{-(n-\nu)/q)} (\lambda_{q+1} \ldots \lambda_n)^{1/q}. \]
We infer 
\[ \gamma_\epsilon := \frac{q\epsilon}{\hat{\lambda}_1 + \ldots + 
                                       \hat{\lambda}_q} \leq 
   \min(1,\frac{q\epsilon}{\delta \lambda_q}) \leq
   \min(1, C^\prime \delta^{-1} \epsilon^{-(n-\nu)/q)} 
                              (\lambda_{q+1} \ldots \lambda_n)^{1/q}). \]
We notice that
\[ \int_X \lambda_{q+1} \ldots \lambda_n dV_\omega \leq 
   \int_X (\Theta_{h_\epsilon}(L) + \epsilon\omega)^{n-q} \wedge \omega^q \leq
   (c_1(L) + \epsilon\{\omega\})^{n-q} \{\omega\}^q \leq C^{\prime\prime}, \]
hence the functions 
$(\lambda_{q+1} \ldots \lambda_n)^{1/q}$
are uniformly bounded in
$L^1$
norm as
$\epsilon$
tends to
$0$. 
Since
$1 - (n-\nu)/q > 0$
by hypothesis, we conclude that 
$\gamma_\epsilon$
converges almost everywhere to
$0$
as
$\epsilon$
tends to 
$0$. 
On the other hand,
\[ |\beta|^2_{\hat{h}_\epsilon} = |\beta|^2_{h_\infty} 
   e^{-2(\delta(1+s)\phi_\epsilon + (1-\delta)(1+s)\psi_\epsilon)} \leq
   |\beta|^2_{h_\infty} e^{-2(\delta(1+s)\phi_\epsilon)}
                        e^{-2(1-\delta)(1+s)\psi}. \]
Our assumption that the coefficients of
$\beta$
lie in
$\mathcal{J}(h^{1+s^\prime}_{\min})$
implies that there exists a
$p^\prime > 1$ 
such that
$\int_X |\beta|^2_{h_\infty} e^{-2p^\prime(1-\delta)(1+s)\psi} dV_\omega 
 < \infty$,
no matter how small
$\delta$
is. Let
$p \in (1,\infty)$
be the conjugate exponent such that
$\frac{1}{p} + \frac{1}{p^\prime} = 1$. 
By H\"older's inequality, we have
\[ \int_X \gamma_\epsilon |\beta|^2_{\hat{h}_\epsilon} dV_\omega \leq
   (\int_X |\beta|^2_{h_\infty} e^{-2p\delta(1+s)\phi_\epsilon} 
                                                       dV_\omega)^{1/p}
   (\int_X \gamma_\epsilon^{p^\prime} |\beta|^2_{h_\infty} 
         e^{-2p^\prime(1+s)(1-\delta)\psi} dV_\omega)^{1/p^\prime}.\]
As
$\gamma_\epsilon \leq 1$,
the Lebesgue dominated convergence theorem shows that
\[ \int_X \gamma_\epsilon^{p^\prime} |\beta|^2_{h_\infty} 
         e^{-2p^\prime(1+s)(1-\delta)\psi} dV_\omega \]
converges to
$0$
as
$\epsilon$
tends to
$0$.

\noindent
For the first integral, we argue as follows: The
$\phi_\epsilon$
may be constructed such that the Lelong numbers 
$\nu((1+s)\phi_\epsilon,x)$
are bounded from above by the Lelong numbers
$\nu((1+s)\psi,x)$
in every point
$x \in X$
(see again Theorem 3.16 in \cite{Eck03}). On the other hand, there is a 
constant
$C$
such that
$\nu((1+s)\psi,x) < C$
for all points
$x \in X$,
by \cite[Lem.3.11]{Bou02}. Hence,
$\nu(\frac{1+s}{C}\phi_\epsilon,x) < 1$,
and
\[ \int_X e^{-(2/C)(1+s)\phi_\epsilon} dV_\omega < \infty, \]
by Skoda's lemma \cite[(5.6)]{Dem00}. Adding sufficiently big constants to
$\phi_\epsilon$
we can even reach that the integrals above are \textit{uniformly bounded}. By
choosing
$\delta \leq 1/(pC)$, 
the integral
$\int_X |\beta|^2_{h_\infty} e^{-2p\delta(1+s)\phi_\epsilon} dV_\omega$
remains bounded and we are done.

\newcommand{\etalchar}[1]{$^{#1}$}


\begin{thebibliography}{BCE{\etalchar{+}}00}

\bibitem[Bha01]{Bha01}
R.~Bhatia.
\newblock Linear algebra to Quantum Cohomology: The Story of Alfred Horn's
  Inequalities.
\newblock {\em Math. Monthly}, 108:289--318, Apr 2001.

\bibitem[Bon95]{Bon95}
L.~Bonavero.
\newblock {In\'egalit\'es de Morse et vari\'et\'es de Moishezon}.
\newblock Preprint, math.AG/9512013, 1995.

\bibitem[Bou02]{Bou02}
S.~Boucksom.
\newblock {\em {C\^{o}nes positifs des vari\'{e}t\'{e}s complexes compactes}}.
\newblock PhD thesis, Grenoble, 2002.

\bibitem[Dem82]{Dem82}
J.-P. Demailly.
\newblock {Estimations $L^2$ pour l'op\'{e}rateur $\overline{\partial}$ d'un
  fibr\'{e} vectoriel holomorphe semi-positif au dessus d'une vari\'{e}t\'{e}
  k\"ahlerienne compl\`{e}te}.
\newblock {\em Ann.Sci. ENS}, 15:457--511, 1982.

\bibitem[Dem00]{Dem00}
J.-P. Demailly.
\newblock Multiplier ideal sheaves and analytic methods in algebraic geometry.
\newblock Lecture Notes, School on Vanishing theorems and effective results in
  Algebraic Geometry, ICTP Trieste, April 2000.

\bibitem[DP02]{DP02}
J.-P. Demailly and Th. Peternell.
\newblock {A Kawamata-Viehweg Vanishing Theorem on compact K{\"a}hler
  manifolds}.
\newblock math:AG/0208021, August 2002.

\bibitem[DPS94]{DPS94}
J.-P. Demailly, Th. Peternell, and M.~Schneider.
\newblock Compact complex manifolds with numerically effective tangent bundles.
\newblock {\em J. Alg. Geom.}, 3:295--345, 1994.

\bibitem[DPS01]{DPS00}
J.-P. Demailly, Th. Peternell, and M.~Schneider.
\newblock Pseudo-effective line bundles on compact k\"ahler manifolds.
\newblock {\em Int. J. Math..}, 12(6):689--741, 2001.

\bibitem[Eck03]{Eck03}
Thomas Eckl.
\newblock {Numerical Trivial Foliations}.
\newblock math.AG/0304312, 2003.

\bibitem[H{\"o}r65]{Hoer65}
L.~H{\"o}rmander.
\newblock {$L^2$ estimates and Existence Theorems for the
  $\overline{\partial}$-\!operator}.
\newblock {\em Acta Math.}, 113:89--152, 1965.

\bibitem[PH78]{GH}
P.Griffiths and J.~Harris.
\newblock {\em Principles of Algebraic Geometry}.
\newblock Wiley, New York, 1978.

\bibitem[SN67]{Nag67}
B.~Sz{\"o}kefalvi-Nagy.
\newblock {\em Spektraldarstellung linearer Transformationen des Hilbertschen
  Raums}.
\newblock Springer, Berlin, 1967.

\bibitem[Wel80]{Wel80}
R.~Wells.
\newblock {\em Differential Analysis on Complex manifolds}.
\newblock Springer, New York, 1980.

\end{thebibliography}

\end{document}